\documentclass[11pt]{article}
\usepackage{amssymb}
\usepackage{amsmath}
\usepackage{amsthm}
\usepackage{latexsym}
\usepackage{amsfonts}
\usepackage{graphicx}
\usepackage{graphics}

\newcommand{\lyxaddress}[1]{
\par {\raggedright #1
\vspace{1.4em}
\noindent\par}
}

\newcommand{\dis}{\displaystyle}
\theoremstyle{plain}
\newtheorem{thm}{Theorem}[section]   
\newtheorem{prop}[thm]{Proposition}
\newtheorem{lem}[thm]{Lemma}
\newtheorem{cor}[thm]{Corollary}

\theoremstyle{definition}
\newtheorem{rem}[thm]{Remark}
\newtheorem{exm}[thm]{Example}
\newtheorem*{Proof}{Proof}

\newcommand{\ra}{\;\rightarrow\;}

\newcommand{\al}{\alpha}

\newcommand{\Ga} {{\varGamma}}

\newcommand{\de}{\delta }
\newcommand{\OO} {{\varOmega}}

\newcommand{\e}{\varepsilon }

\newcommand{\C}{\mathbb{C}}
\newcommand{\R}{\mathbb{R}}

\newcommand{\N}{\mathbb{N}}

\newcommand{\tf}{\widetilde{f}}
\newcommand{\tg}{\widetilde{g}}

\newcommand{\ld}{\ldots}

\newcommand{\sm}{\smallsetminus}

\newcommand{\qb}{$\quad\blacksquare$}
\begin{document}
\title{\bf On a characterization of Arakelian sets}
\author{G. Fournodavlos}

\maketitle

\lyxaddress{Department of Mathematics, University of Athens, Panepistemioupolis,
15784, Athens, Greece, e-mail: gregdavlos@hotmail.com}
%
\noindent
%
\begin{abstract}
Let $K$ be a compact set in the complex plane $\C$, such that its complement in the Riemann sphere,
$(\C\cup\{\infty\})\sm K$, is connected. Also, let $U\subseteq\C$ be an open set which contains $K$.
Then there exists a simply connected open set $V$ such that $K\subseteq V\subseteq U$.
We show that if the set $K$ is replaced by a closed set $F$ in $\C$, then the above lemma is equivalent to the fact that $F$
is an Arakelian set in $\C$. This holds more generally, if $\C$ is replaced by any simply connected open set
$\OO\subseteq\C$. In the case of an arbitrary open set $\OO\subseteq\C$, the above extends to the one point
compactification of $\OO$. As an application we give a simple proof of the fact that the disjoint union of two
Arakelian sets in a simply connected open set $\OO$ is also Arakelian in $\OO$.
\end{abstract}
{\em Subject Classification MSC2010}\,: 30E10 30K05 \vspace*{0.2cm} 
\section{Introduction}\label{sec1} 
\noindent

K.-G. Grosse-Erdmann in \cite{13} and G. Costakis in \cite{4} stated and proved, independently, the following.
\begin{lem}\label{lem1.1}
Let $K\subseteq\C$ be a compact set with $(\C\cup\{\infty\})\sm K$ connected. If $U$ is an open set in
$\C$ containing $K$, then there exists a simply connected open set $V$ such that $K\subseteq V\subseteq U$.
\end{lem}

The above lemma gave several applications. In particular we used it recently in \cite{7}. I was interested in
knowing if the above lemma still holds if the compact set $K$ is replaced by a closed subset $F$ in $\C$.
The answer is negative in general. A counterexample is the following:
\[
F=\dis\bigcup^\infty_{n=1}\bigg[\big(\{\dis\sum^n_{i=1}\frac{1}{2^i}\}\times[0,n]\big)
\cup\big([\dis\sum^n_{i=1}\frac{1}{2^i},\dis\sum^{n+1}_{i=1}\frac{1}{2^i}]\times\{n\}\big)\bigg]
\cup\big(\{1\}\times[0,\infty)\big).
\]
This set $F$ relates to the well-known Arakelian's Approximation Theorem \cite{1}.
\begin{thm}\label{thm1.2}
Let $F$ be a closed set in the complex plane $\C$. Then every function $f:F\ra\C$ continuous on $F$
and holomorphic in $F^0$ ($f\in A(F)$) can be uniformly approximated on $F$ by entire functions, $g\in H(\C)$,
if and only if the following hold:

(i)$(\C\cup\{\infty\})\sm F$ is connected

(ii)$(\C\cup\{\infty\})\sm F$ is locally connected at $\infty$.
\end{thm}

Yet, in \cite{14} one finds another proof of Theorem \ref{thm1.2}, based on Mergelyan's theorem,
where conditions $(i)$ and $(ii)$ are replaced by the next equivalent condition:

$(iii)$ $(\C\cup\{\infty\})\sm F$ is connected and for every closed disk $D$, in $\C$, the union
of all bounded components of $\C\sm(F\cup D)$ is bounded.
A set satisfying $(iii)$ (or equivalently $(i)$ and $(ii)$) is called an Arakelian set in $\C$.
In the present article we prove the following.
\begin{thm}\label{thm1.3}
Let $F$ be a closed subset of $\C$. Then the following are equivalent:

(1) For every open set $U\subseteq\C$, which contains $F$, there exists a simply connected open
set $V$ such that $F\subseteq V\subseteq U$

(2)$F$ is an Arakelian set in $\C$.
\end{thm}

More generally, Arakelian sets may be defined for the arbitrary open set $\OO\subseteq\C$.
The question is for a relatively closed set $F$ in $\OO$ whether every function $f\in A(F)$
can be uniformly approximated on $F$ by holomorphic functions $g\in H(\OO)$. This was
completely settled by Arakelian in \cite{2}, where he extended Theorem \ref{thm1.2}. In this
version one considers the one point compactification of $\OO$, $\OO\cup\{\al\}$. The relatively
closed set $F\subseteq\OO$ for which the approximation is possible is called Arakelian set in $\OO$
and again a purely topological description is possible.

We extend Theorem \ref{thm1.3} replacing the complex plane $\C$ by any simply connected open set
$\OO\subseteq\C$. This means that the open set $V$, in the extended version of Theorem \ref{thm1.3} (1),
is still simply connected; that is $(\C\cup\{\infty\})\sm V$ is connected. Also, we can further extend Theorem
\ref{thm1.3} in the general case of any open set $\OO\subseteq\C$, not necessarily simply connected.
In this case the open set $V$ is not simply connected, but its complement $(\OO\cup\{\al\})\sm V$ in the
one point compactification of $\OO$ has to be connected.

Next, we give two applications of our results. One of these states that if $\OO$ is a simply connected open
subset of $\C$, then the union of two disjoint Arakelian sets in $\OO$ is also Arakelian in $\OO$. We notice
that when $\OO$ is not simply connected the above fails.

N. Tsirivas in \cite{16} proved a variation of Lemma \ref{lem1.1}, without the assumption that
$(\C\cup\{\infty\})\sm K$ is connected; see also \cite{6}. There are indications that we could obtain similar variations
corresponding to our results, but we have not yet managed to do so.

Some partial extensions of Theorem \ref{thm1.3} can be achieved in the case of Riemann surfaces (\cite{8}).
This is another possible direction for further investigation in connection with our results.
P. M. Gauthier suggested that some alternative proofs could relate to Runge's pairs
and harmonic approximation \cite{3}, \cite{9}, \cite{10}, \cite{11}.

Finally, we mention that it is open to characterize all subsets $F\subseteq\C$, such that the conclusion
of Lemma \ref{lem1.1} holds. We can easily find examples of sets $F$, which are not relatively closed in
any simply connected open set $\OO\subseteq\C$; in some of these examples the conclusion of Lemma \ref{lem1.1}
holds, but in some others it does not.
\section{The results}\label{sec2}
\noindent

In \cite{2} N. U. Arakelian proved the following theorem.

\begin{thm}\label{thm2.1}
Let $\OO\subseteq\C$ be an open set and $F$ a relatively closed subset of $\OO$. Then
every function $f\in A(F)$ can be uniformly approximated on $F$ by functions $g\in H(\OO)$
holomorphic in $\OO$, if and only if the following hold:

i) $(\OO\cup\{\al\})\sm F$ is connected and

ii) $(\OO\cup\{\al\})\sm F$ is locally connected at $\al$, where $\OO\cup\{\al\}$ is the one
point compactification of $\OO$.

Such a set $F$ is also called an Arakelian set in $\OO$.
\end{thm}

We say that $B$ is a "hole of $F$" in $\OO$, iff $B$ is a component of $\OO\sm F$,
which is contained in a compact subset of $\OO$. Note that
$(\OO\cup\{\al\})\sm F$ is connected iff $F$ has no holes in $\OO$.
\begin{prop}\label{prop2.2}
A closed set $F$ in $\OO$, without holes, is an Arakelian set in $\OO$, if and only if
for every compact set $K\subseteq\OO$, the union of all holes of $F\cup K$ in $\OO$,
is contained in a compact subset of $\OO$.
\end{prop}

For the case $\OO=\C$ see \cite{14}. We include a proof of the general case, for the sake of completeness.
\begin{Proof}$(\Rightarrow)$
Suppose that there is a compact set $K\subseteq\OO$, such that the union of all holes in
$\OO$ of $F\cup K$ is either unbounded or has zero distance from the boundary of $\OO$. The complement
$(\OO\cup\{\al\})\sm K$ is a neighborhood of $\al$, in $\OO\cup\{\al\}$. Hence, there exists a neighborhood
$W\subseteq(\OO\cup\{\al\})\sm K$ of $\al$, in $\OO\cup\{\al\}$, such that
$W\cap\big[(\OO\cup\{\al\})\sm F\big]$ is connected. There
exists a hole $B$ of $F\cup K$ such that $B\cap W\neq\emptyset$, because $(\OO\cup\{\al\})\sm W$ is contained in
a compact subset of $\OO$. Observe that $bd(B)\subseteq F\cup K$. Since $W\cap\big[(\OO\cup\{\al\})\sm F\big]$
and $F\cup K$ are disjoint, it follows that
$W\cap\big[(\OO\cup\{\al\})\sm F\big]\subseteq B\cup\big[(\OO\cup\{\al\})\sm\overline{B}\big]$, which is a contradiction.

$(\Leftarrow)$ Let $U\subseteq\OO\cup\{\al\}$ be an open neighborhood of $\al$ in $\OO\cup\{\al\}$.
The set $K=(\OO\cup\{\al\})\sm U$ is compact. Therefore, the union of all holes in $\OO$ of $F\cup K$
is contained in a compact subset of $\OO$.
Let $B_1,B_2\ld$ be those holes. Also, let $W=(\OO\cup\{\al\})\sm(K\cup B_1\cup B_2\cup\cdots)$. Obviously, $W$ is
a neighborhood of $\al$ and $W\subseteq U$. We notice that $W\cap\big[(\OO\cup\{\al\})\sm F\big]$ is the union
of $\{\al\}$ and all the components of $\OO\sm(F\cup K)$, which are either unbounded or have zero distance from
$bd(\OO)$. Thus, $W\cap\big[(\OO\cup\{\al\})\sm F\big]$ is connected and the proof is complete.\qb
\end{Proof}
\begin{rem}\label{rem2.3}
In order to determine whether a relatively closed set $F$, without holes, is Arakelian in $\OO$,
it suffices to check the condition of Proposition \ref{prop2.2} only
for an exhausting sequence, $(K_n)_{n\in\N}$, of compact subsets of $\OO$. Such a sequence can be chosen
so that $K_n\subseteq K^0_{n+1}$, $n\in\N$, $\dis\bigcup^\infty_{n=1}K^0_n=\OO$
and $K_n$ has no holes in $\OO$, for all $n\in\N$ (\cite{15}). Also, we can assume that $K_n$ is a finite
union of squares in a grid, whose sides are parallel to the coordinate axes
and of length $\de_n>0$, $n\in\N$.
\end{rem}

\begin{thm}\label{thm2.4}
If $F$ is an Arakelian set in $\OO$, then for every open set $U\subseteq\OO$, which contains $F$, there exists
an open set $V\subseteq\OO$ such that $F\subseteq\ V\subseteq U$ and $(\OO\cup\{\al\})\sm V$ is connected.
\end{thm}
\begin{Proof}
Let $F$ be an Arakelian set in $\OO$ and $U\subseteq\OO$ an open set such that $F\subseteq U$.
We define $d_x=\min\{\frac{dist(x,F)}{2}, dist(x,\C\sm\OO)), 1\}>0$, $x\in\OO\sm U$. Also, let $(K_n)_{n\in\N}$ be an
exhausting sequence of compact subsets of $\OO$.
\begin{itemize}
\item The relatively closed set $\OO\sm U$ has a locally finite cover in $\OO$,
$\big\{D(x_i,d_{x_i})\big\}^\infty_{i=1}$. It suffices to choose a finite cover of disks $D(x,d_x)$,
$x\in\OO\sm U$, for each of the compact sets $(K_n\sm K^0_{n-1})\cap(\OO\sm U)$, $n\in\N$, $K_0=\emptyset$.
\end{itemize}

This implies that every compact subset of $\OO$ intersects a finite number of disks from
$\big\{D(x_i,d_{x_i})\big\}^\infty_{i=1}$. Observe that
$\dis\bigcup^\infty_{i=1}\overline{D(x_i,d_{x_i})}$ is closed in $\OO$
and the set $\{x_i$ $|$ $i=1,2,\ld\}$ has no accumulation points in $\OO$. Hence,
$U_1=U\sm\big(\dis\bigcup^\infty_{i=1}\overline{D(x_i,d_{x_i})}\big)$
is open and $F\subseteq U_1\subseteq U$.
\begin{itemize}
\item We say that a point $x\in\OO$ is joined with $\al$ by a curve $\Ga$ in $E\subseteq\OO$, if
$\Ga:[0,+\infty)\ra E$ is continuous and $\Ga(0)=x$, $\dis\lim_{t\to+\infty}\Ga(t)=\al$. The image
of such a curve, $\Ga\big([0,+\infty)\big)$, is relatively closed in $\OO$
\end{itemize}

Each $x_i$ can be joined with $\al$ by a curve $\Ga_i$, $i=1,2,\ld$, such that for every $n\in\N$ only finite
curves intersect the compact set $K_n$. Indeed, if $n\in\N$, then there are finitely many $x_i$ contained in the union
of $K_n$ and all the holes of $F\cup K_n$ in $\OO$. The points that we have not already joined
with $\al$ (induction), are contained in components of $\OO\sm(F\cup K_{n-1})$, which are either
unbounded or have zero distance from the boundary of $\OO$. Let $x_i$ be such a point and $E$ the component
of $\OO\sm(F\cup K_{n-1})$, which contains it. We can construct a curve $\Ga_i$ in $E$, which joins $x_i$ with $\al$.
\begin{itemize}
\item If $s\ge n$, then by Proposition \ref{prop2.2} $E$ contains a component
$E_s$ of $\OO\sm(F\cup K_n)$, which is either unbounded or has zero distance from $bd(\OO)$.
Further, we can assume that $E_{s-1}\supseteq E_s$, $s\ge n$, where $E_{n-1}=E$.
Let $x_{is}\in E_s$ and $\Ga_{is}$ a curve in $E_{s-1}$, which joins $x_{is-1}$
with $x_{is}$ (such a curve exists, since $E_{s-1}$ is open and connected), $s\ge n$,
where $x_{in-1}=x_i$. The desired curve, $\Ga_i$, consists of all $\Ga_{is}$, $s\ge n$.
\end{itemize}

Thus, the union $\dis\bigcup^\infty_{i=1}\Ga_i$ is closed in $\OO$ and the open set
$V=U_1\sm(\dis\bigcup^\infty_{i=1}\Ga_i)$ has the desired properties. Obviously,
$F\subseteq V\subseteq U_1\subseteq U$ and
$(\OO\cup\{\al\})\sm V=\dis\bigcup^\infty_{i=0}\big(\overline{D(x_i,d_{x_i})}\cup\Ga_i\big)\cup\{\al\}$
is connected, which completes the proof.
\qb
\end{Proof}

We note that the previous theorem, in the case $\OO=\C$, is known; see \cite{8}, \cite{12}.
\begin{thm}\label{thm2.5}
If $F$ is a closed set in $\OO$, such that for every open set $U\subseteq\OO$, which contains $F$,
there exists an open set $V\subseteq\OO$ with $F\subseteq V\subseteq\OO$ and $(\OO\cup\{\al\})\sm V$ connected,
then $F$ is an Arakelian set in $\OO$.
\end{thm}
\begin{Proof}
First, we notice that $F$ has no holes in $\OO$. If $B$ is a hole of $F$ in $\OO$ and $x\in B$, then the open set
$U=\OO\sm\{x\}$ contains $F$. Hence, there exists an open set $V\subseteq\OO$ such that $F\subseteq V\subseteq U$ and
$(\OO\cup\{\al\})\sm V$ is connected. It holds $(\OO\cup\{\al\})\sm F\supseteq(\OO\cup\{\al\})\sm V$.
Therefore, the latter is contained in the component of $(\OO\cup\{\al\})\sm F$ that contains $\al$. However,
$\big[(\OO\cup\{\al\})\sm V\big]\cup B\neq\emptyset$, which is a contradiction, because $B$ is a component of
$(\OO\cup\{\al\})\sm F$ not containing $\al$.

Suppose that $F$ is not an Arakelian set in $\OO$. By Proposition \ref{prop2.2} there exists
a compact set $K\subseteq\OO$, such that the union of all holes of $F\cup K$ in $\OO$, is either unbounded or has
zero distance from $bd(\OO)$. Moreover, Remark \ref{rem2.3} enables us to assume that $K$ is a finite
union of closed squares in a grid, whose sides are parallel to the coordinate axes and of length $\de>0$.
Let $B_1,B_2,\ld$ be a sequence of holes of $F\cup K$ in $\OO$ and
$x_n\in B_n$, $n\in\N$, such that $x_n\ra\al$, as $n\ra+\infty$.

The open set $U=\OO\sm\{x_1,x_2,\ld\}$ contains $F$. Thus, there
exists an open set $V$ with $F\subseteq V\subseteq U$ and $(\OO\cup\{\al\})\sm V$ connected.
Observe that $(\OO\cup\{\al\})\sm V$ intersects $B_n$ and $(\OO\cup\{\al\})\sm\overline{B_n}$, for all $n\in\N$.
This implies that there exists $y_n\in(\OO\sm V)\cap bd(B_n)\cap bd(K)$, $n\in\N$.
Since $bd(K)$ is compact, $(y_n)_{n\in\N}$ has a limit point $y\in bd(K)$. Also, $y\in\OO\sm V$, because
$\OO\sm V$ is closed in $\OO$ and $y\in\OO$. We claim that $y\in F\subseteq V$, which is obviously
a contradiction. Indeed, if $y\not\in F$, then there exists $\e>0$ such that
$D(y,\e)\subseteq\OO$ does not intersect $F$. In addition, we can choose $\e>0$,
depending on the place of $y$ in the grid, so that $D(y,\e)\sm K$ has at most two components.
This is a contradiction, since $D(y,\e)\sm K\subseteq\OO\sm(F\cup K)$ intersects infinite holes from
$\big\{B_n\big\}^\infty_{n=1}$. The proof is complete.\qb
\end{Proof}

According to Theorems \ref{thm2.5} and \ref{thm2.6}, we have the following characterization
of Arakelian sets.
\begin{thm}\label{thm2.6}
Let $\OO\subseteq\C$ be an open set and $\OO\cup\{\al\}$ its one point compactification.
A relatively closed set $F$ is Arakelian in $\OO$, if and only if
for every open set $U\subseteq\OO$, which contains $F$, there exists an open set $V\subseteq\OO$ such that
$F\subseteq V\subseteq U$ and $(\OO\cup\{\al\})\sm V$ is connected.
\end{thm}
\begin{lem}\label{lem2.7}
Let $\OO\subseteq\C$ be a simply connected open set. A set $G\subseteq\OO$ has
connected complement in $\OO\cup\{\al\}$, if and only if its complement in the Riemann sphere,
$(\C\cup\{\infty\})\sm G$, is connected.
\end{lem}
\begin{Proof}$(\Rightarrow)$
Let $G\subseteq\OO$ with $(\OO\cup\{\al\})\sm G$ connected. Assume that $(\C\cup\{\infty\})\sm G$ is
not connected. Thus, there are two open sets $U_1,U_2$ in $\C\cup\{\infty\}$, such that
$U_i\cap\big[(\C\cup\{\infty\})\sm G\big]\neq\emptyset$, $i=1,2$,
$(\C\cup\{\infty\})\sm G\subseteq U_1\cup U_2$ and $U_1\cap U_2\cap\big[(\C\cup\{\infty\})\sm G\big]=\emptyset$.
Since $(\C\cup\{\infty\})\sm\OO$ is connected and it is contained in $(\C\cup\{\infty\})\sm G$,
it follows that $(\C\cup\{\infty\})\sm\OO$ is contained in exactly one of the sets $U_1,U_2$.
Without loss of generality, we assume that $(\C\cup\{\infty\})\sm\OO\subseteq U_1$.
Observe that $(\C\cup\{\infty\})\sm U_1$ is a compact subset of $\OO$. This implies
that $V_1=(U_1\cap\OO\sm G)\cup\{\al\}$ and the set
$V_2=U_2\cap\big[(\OO\cup\{\al\})\sm G\big]\subseteq(\C\cup\{\infty\})\sm U_1\subseteq\OO$ are two
nonempty disjoint open sets in $(\OO\cup\{\al\})\sm G$.
Furthermore, it is easy to see that $(\OO\cup\{\al\})\sm G\subseteq V_1\cup V_2$ and thus we obtain a contradiction.

$(\Leftarrow)$ Let $G\subseteq\OO$ with $(\C\cup\{\infty\})\sm G$ connected. We define
$\phi:\C\cup\{\infty\}\ra\OO\cup\{\al\}$,
$\phi(x)=\left\{
  \begin{array}{ll}
    x, & x\in\OO \\
    \al, & x\not\in\OO.
  \end{array}
\right.$ Obviously, $\phi$ is continuous and so
$\phi\big((\C\cup\{\infty\})\sm G\big)=(\OO\cup\{\al\})\sm G$ is connected.
\qb
\end{Proof}

Combining Theorem \ref{thm2.6} with Lemma \ref{lem2.7}, we obtain the following.
\begin{thm}\label{thm2.8}
Let $\OO\subseteq\C$ be a simply connected open set and $F\subseteq\OO$ a relatively closed set. 
Then the following are equivalent:

i)$F$ is an Arakelian set in $\OO$

ii)For every open set $U\subseteq\C$, which contains $F$, there exists a simply connected open set
$V\subseteq\C$ such that $F\subseteq V\subseteq U$.
\end{thm}

The next corollary is an immediate application of Theorem \ref{thm2.8}.
\begin{cor}\label{cor2.9}
If $\OO\subseteq\C$ is a simply connected open set, then the disjoint union of
two Arakelian sets in $\OO$ is also Arakelian in $\OO$.
\end{cor}
\begin{Proof}
Let $F_1,F_2$ be two disjoint Arakelian sets in $\OO$. Also, let $U\subseteq\OO$ an open set,
which contains the union $F_1\cup F_2$. Since $F_1$ and $F_2$ are two
disjoint closed sets in $\OO$, there exist two disjoint open sets $G_1,G_2\subseteq\OO$ such that
$F_1\subseteq G_1$ and $F_2\subseteq G_2$. By Theorem \ref{thm2.8} there two simply connected
open sets $V_1,V_2$ with $F_i\subseteq V_i\subseteq G_i\cap U$, $i=1,2$. Obviously, it holds
$F_1\cup F_2\subseteq\ V_1\cup V_2\subseteq U$ and since $V_1\cap V_2=\emptyset$, every
component of $V_1\cup V_2$ is simply connected. This implies that $V_1\cup V_2$ is a simply
connected open set. Thus, according to theorem Theorem \ref{thm2.8},
the closed set $F_1\cup F_2$ is Arakelian in $\OO$ and the proof is complete.
\qb
\end{Proof}

In the case $\OO=\C$, an alternative proof of the previous result, using Proposition \ref{prop2.2},
can be found in \cite{5}. The following example shows that Corollary \ref{cor2.9} does not hold when $\OO$ is
not simply connected.
\begin{exm}\label{exm2.10}
Let $\OO=D(0,1)\sm\{0\}$. Also, let $F_1=C(0,r_1)$ and $F_2=C(0,r_2)$, where $0<r_1<r_2<1$.
Observe that $F_1,F_2$ are two disjoint compact subsets of $\OO$ with connected complements in the
one point compactification of $\OO$. Hence, both sets are Arakelian in $\OO$. Nonetheless, the union $F_1\cup F_2$
is not Arakelian in $\OO$, since $(\OO\cup\{\al\})\sm(F_1\cup F_2)$ is not connected.
\qb
\end{exm}

Even if $\OO\subseteq\C$ is a simply connected open set, it is not true that
the infinite denumerable union of pairwise disjoint Arakelian sets in $\OO$ is also Arakelian in $\OO$.
\begin{exm}\label{exm2.11}
Let $\OO=\C$ and let $F_0=\{2\}\times\R$, $F_n=\big(\{\dis\sum^{n-1}_{i=0}\frac{1}{2^i}-\frac{1}{2^n},\dis\sum^{n-1}_{i=0}\frac{1}{2^i}\}\times[0,n]\big)
\cup\big([\dis\sum^{n-1}_{i=0}\frac{1}{2^i}-\frac{1}{2^n},\dis\sum^{n-1}_{i=0}\frac{1}{2^i}]\times\{n\}\big)$, $n\ge1$.
It is easy to see that each $F_n$, $n=0,1,\ld$, is an Arakelian set in $\C$. However, $F=\dis\bigcup^\infty_{n=0}F_n$
is not Arakelian in $\C$, because despite the fact that $F$ is closed and $(\C\cup\{\infty\})\sm F$ is connected, the union
of all holes in $\C$ of $\overline{D(0,r)}\cup F$, $r\ge2$, is unbounded.
\qb
\end{exm}

Finally, we present another application of our characterization.
\begin{cor}\label{cor2.12}
Let $\OO\subseteq\C$ be a simply connected open set and $F\subseteq\OO$ a relatively closed set. Also, let
$f\in A(F)$ satisfying $f(z)\neq0$, for all $z\in F$. Then there exists a function $g\in A(F)$ such that
$f=e^g$.
\end{cor}
\begin{Proof}
According to Tietze's extension theorem, there exists a continuous extension of $f$ on $\OO$, which we
denote by $\tf:\OO\ra\C$. The open set $U=\OO\sm\tf^{-1}(0)$ contains $F$. By Theorem \ref{thm2.8},
there is a simply connected open set $V$ with $F\subseteq\ V\subseteq U$.
This implies that there exists a continuous function $\tg:V\ra\C$ such that
$\tf_{\big|V}=e^{\tg}$. The function $g=\tg_{\big|F}$ is obviously continuous on $F$ and $f=e^g$.
Since $f_{\big|F^0}$ is holomorphic, $g$ is also holomorphic in $F^0$. Thus, $g\in A(F)$ and the proof is complete.
\qb
\end{Proof}

We notice that for $\OO=\C$ Corollary \ref{cor2.12} is known; see \cite{12}.

\bigskip

{\bf Acknowledgement}: I would like to express my thanks to V. Nestoridis
for his valuable suggestions and his interest in this work. Also, I would like
to thank P. M. Gauthier for his interest in this work, for bringing to my attention
references I was not aware of and for suggestions towards further work.

\end{document}